\documentclass{article}

\usepackage[pdftex]{graphicx}
\usepackage{amsmath,amssymb,fullpage,enumerate,amsmath,amsthm}

\numberwithin{equation}{section}
\usepackage{tikz}\usetikzlibrary{arrows}
\usepackage{subfigure}

\usepackage{authblk} 

\usepackage{cite, mdframed}

\newcommand{\av}[1]{\left|{#1}\right|}

\newtheorem*{theorem*}{Theorem}
\newtheorem*{cl*}{Claim}

\theoremstyle{definition}

\theoremstyle{definition}

\theoremstyle{remark}

\theoremstyle{definition}

\newcommand{\E}{\mathbb{E}}
\newcommand{\Var}{{\rm{Var}}}

\newcommand{\mcL}{\mathcal{L}}

\title{Finite-size effects and switching times for Moran dynamics with
  mutation} \author{Lee DeVille and Meghan Galiardi} \affil{Department
  of Mathematics, University of Illinois}

\begin{document}
\maketitle

\begin{abstract}
  We consider the Moran process with two populations competing under
  an iterated Prisoners' Dilemma in the presence of mutation, and
  concentrate on the case where there are multiple Evolutionarily
  Stable Strategies.  We perform a complete bifurcation analysis of
  the deterministic system which arises in the infinite population
  size.  We also study the Master equation and obtain asymptotics for
  the invariant distribution and metastable switching times for the
  stochastic process in the case of large but finite population.  We
  also show that the stochastic system has asymmetries in the form of
  a skew for parameter values where the deterministic limit is
  symmetric.

\end{abstract}

\section{Introduction}

The mathematical study of models for evolution has a long history.
The early work of Wright, Fisher, Moran, and Feller~\cite{Wright.31,
  Feller.51, Fisher.book, Moran.58}, where those authors studied of
small population size on a model of neutral evolution
(see~\cite{Crow.10} for a nice overview of the early history of these
studies).  In parallel, the mathematical framework of game
theory~\cite{von1953theory} and evolutionary game
theory~\cite{smith1982evolution} has been used since to model and
understand complexity in problems in evolution and ecology.  A truly
large volume of work has been developed in the last few decades using
this approach~\cite{may1975nonlinear, axelrod1981evolution,
  Axelrod.book, Nowak.92, Nowak.book, axelrod1997complexity, Nowak.06,
  Nowak.Sigmund.92, Nowak.Sigmund.05, Taylor.04,
  weibull1997evolutionary}.

Most of the models mentioned above can be formulated as a Markov chain
on a lattice, where the discrete points on the lattice correspond to
integer sizes of the various subpopulations.  In the limit of large
population, this Markov chain undergoes smaller jumps (on the
population scale), if this is coupled with an appropriate rescaling in
time, we might expect a sensible limit to a continuum model.  What is
true is that if one does this rescaling properly, then in the limit of
infinite population, the system does converge to an ODE.  In the case
of large but finite populations, the system is still stochastic, so it
is of interest to study the invariant distributions of such a system,
the rate and type of convergence to this equilibrium, and in
particular switching processes between multiple equilibria.  The first
rigorous results in this direction were due to the pioneering work of
Kurtz in the context of biochemical models~\cite{Kurtz.72, Kurtz.78},
but have been followed up by a large number of analytic and numerical
studies~\cite{Bressloff.book, Bressloff.10, vanKampen.book,
  vanKampen.82, Darling.Norris.08, Gardiner.book, SW,
  guttenberg2008cascade}.

The purpose of this paper is to apply the modern techniques of
stochastic analysis to a particular model of evolutionary game theory.
The model we consider is a version of a model described
in~\cite{Nowak.book}; we consider a population consisting of
organisms who can each play one of two strategies where the payoff
matrix is of Prisoner's Dilemma type, and the organisms evolve
according to Moran dynamics, i.e. in each round of play, one organism
is chosen to reproduce based on its fitness and one organism is killed
at random, so that the population size remains constant.  We also
allow a positive probability for a mutation, i.e. there is a positive
probability $\mu>0$ that an organism's offspring is of the other type.
We stress here that we consider a fully stochastic version of this
model where all events are governed by a Markov process.

The infinite-population limit of this process (an ODE) is described
in~\cite{Nowak.book}, and it is shown there that for certain
parameters the ODE supports two stable fixed points with disjoint
basins of attraction, i.e. is multistable.  We show that our Markov
chain model limits on this ODE model as the population tends to
infinity.  

For population large but finite, however, the system is still noisy
and can jump between the multiple attractors, although it will only do
so on long timescales.  Studying the timescale on which this jumping
occurs is one of the main goals of this paper, and we give two
different asymptotic descriptions of this timescale.  Another
interesting effect we show by considering the stochastic system (as
compared to its deterministic limit) is we show that certain
parameters that are indistinguishable in the limit give rise to biases
for large but finite populations.  For example, if we assume that the
aggregate payoff is the same for both strategies (i.e. the sum of the
payoffs against both friends and enemies is the same for both
strategies), but one strategy plays better against enemies than it
does against its friends, it will have an advantage for finite size
populations, while being indistinguishable in the limit.

\section{Model} \label{sec:model}
\subsection{Definition} \label{ssec:definition} In this paper, we
consider a finite population of size $N$ consisting of two strategies,
$A$ and $B$.  We assume that payoffs from each individual game are
awarded by a payoff matrix given by
\begin{equation} \label{eq:payoff}
\left[\begin{array}{c|cc}
    & A & B \\\hline
    A & a & b \\
    B & c & d
\end{array}\right].
\end{equation}
Here, we represent the payoff to player1 in an interaction where
player1's play is along the left column, and player2's play is on the
top row.  We assume that the payoff to player2 is always the transpose
of this matrix. For example, when we have two players playing strategy
$A$, they each receive a payoff of $a$, and if they both are strategy
$B$, they each receive $d$.  If two individuals with different
strategies meet, if player1 is type $B$, player1 receives a payoff of
$c$ while player2 receives a payoff of $b$.

Once the average payoffs are determined, each strategy is awarded a
fitness value.  Fitness can have many interpretations such as physical
fitness or ability to collect food, but ultimately it affects the
survival of individuals, and here we assume that the fitness is
proportional to the payoff an individual receives playing the game.
Specifically, if we assume that there are $i$ individuals of type $A$
in the population, then the function
\begin{equation}\label{eq:defoff}
  f_A = \frac{ai + b(N-i)}{N}, \quad\mbox{resp.}\quad f_B = \frac{ci + d(N-i)}{N},
\end{equation}
represents the average fitness of an individual of type $A$
(resp.~type $B$).  This represents the average payoff an organism
would receive if it chose an opponent at random and played the game
described above.

We want to consider only evolutionarily stable strategies (ESS), those
strategies which are resistant to invasion by a single mutant.  For
example, strategy $A$ will be ESS if $f_A>f_B$ when $i=N-1$ (i.e. when
all but one organism is type $A$).  This means that
\begin{equation*}
  a\frac{N-1}N + \frac b N > c\frac{N-1}N + \frac dN,
\end{equation*}
and if we want this equation to hold for arbitrarily large
populations, then the condition is clearly $a>c$.  Similarly, strategy
$B$ is ESS iff $d > b$.  We will assume this in all that follows
below.

After the fitness of each strategy is determined, the population
evolves similar to a Moran Model, but with the possibility of a
mutation during birth. In each round, one individual is chosen for
reproduction with probability proportional to its fitness and gives
birth to one offspring.  However, this offspring is of the opposite
type with probability $\mu\in(0,1)$.  At each stage, one organism is
chosen at random to be removed.

From this, we can view the model as a Markov process $\{X(t)_{t \geq
  0}\}$ on the discrete state space $i = 0, 1, \cdots, N$ with
transition rates
\begin{align}
\omega_+(i;\mu) &= \bigg(\frac{if_A}{if_A + (N-i)f_B}\bigg)\bigg(\frac{N-i}{N}\bigg) (1-\mu) + \bigg(\frac{(N-i)f_B}{if_A + (N-i)f_B}\bigg)\bigg(\frac{N-i}{N}\bigg) \mu \label{eq:transitionRate+} \\
\omega_-(i;\mu) &= \bigg(\frac{(N-i)f_B}{if_A + (N-i)f_B}\bigg)\bigg(\frac{i}{N}\bigg) (1-\mu) + \bigg(\frac{if_A}{if_A + (N-i)f_B}\bigg)\bigg(\frac{i}{N}\bigg) \mu \label{eq:transitionRate-}.
\end{align}
Here the state represents the number of organisms of type $A$.  To see
these rates, note that to have the transition $i \mapsto i+1$, either
an individual with strategy $A$ needs to reproduce (with no mutation)
and an individual with strategy $B$ needs to die, or an individual
with strategy $B$ needs to reproduce (with mutation) and an individual
with strategy $B$ needs to die. Similar arguments can be made for
moving from state $i$ to $i-1$.  (Of course, these rates do not add to
one, since there is a possibility that the population size does not
change: for example, if an individual of type $A$ reproduces without
mutation, but an individual of type $A$ is also selected to die, then
there is no change in the population, etc.)

Throughout this paper we will often switch between the extensive
variable $i$ and the scaled variable $x = \frac{i}{N}$.  In the limit
as $N\rightarrow \infty$, we treat $x$ as a continuous variable. It is
also assumed that the transition rates obey the scaling law
$\omega_{\pm}(i;\mu)=N \Omega_{\pm}(x; \mu)$. 

The following sections study this process both deterministically and
stochastically. For given parameters $a$, $b$, $c$, $d$, and $\mu$ we
seek to describe the distribution of strategies within the population
and determine which strategy (if any) has an advantage.

We also remark that while the definition above seems to average out
some of the randomness before we write down the Markov chain (e.g. we
consider ``average fitness'' when writing down the transition rates in
the chain), we could obtain the same model without such averaging.
For example, if we assume that, at each timestep, each individual is
paired with another individual in the population (selected with
replacement) and then receives a deterministic payoff
from~\eqref{eq:payoff}, and then we choose individuals to reproduce
with a probability proportional to their fitness, then again we would
obtain~(\ref{eq:transitionRate+},~\ref{eq:transitionRate-}).

\subsection{Parameter regimes} \label{ssec:cases}

The results fall into three cases which arise naturally from the payoff matrix \eqref{eq:payoff}. Due to symmetry, we discuss only case 1.1 and case 2.
\begin{itemize}
\item Case 1 --- $a+b = c+d$\\
When populations $A$ and $B$ have equal size, $A$ and $B$ have equal fitness. There are three additional subcases.
\begin{itemize}
\item Case 1.1 --- $a > d$ and $b < c$ \\
$A$ has better self-interactions than strategy $B$, but strategy $B$ has better opponent interactions than strategy $A$.
\item Case 1.2 --- $a < d$ and $b > c$ \\
Same as 1.1, but with the roles of $A$ and $B$ switched
\item Case 1.3 --- $a = d$ and $b = c$ \\
Degenerate case where $A$ and $B$ are interchangable
\end{itemize}
\item Case 2 --- $a+b > c+d$ \\
 When populations $A$ and $B$ have equal sizes,  $A$ has a higher fitness than $B$.
\item Case 3 --- $a+b < c+d$ \\
 When populations $A$ and $B$ have equal sizes,  $B$ has a higher fitness than $A$.
\end{itemize}

For case 2, it is clear that strategy $A$ will have the advantage
based on the asymmetry of the payoff matrix~\eqref{eq:payoff}, and
similarly for Case 3. Case 1.1 is not at all clear {\em a priori}. One
could argue that strategy $A$ has the advantage because it has better
self-interactions than strategy $B$, or that strategy $B$ has the
advantage because it has better interaction with strategy $A$ than
strategy $A$ has with strategy $B$. We will consider these more fully
below.

\subsection{Simulation Results} \label{ssec:simulation}

\begin{figure}[ht]
\begin{centering}
  \includegraphics[width=.48\textwidth]{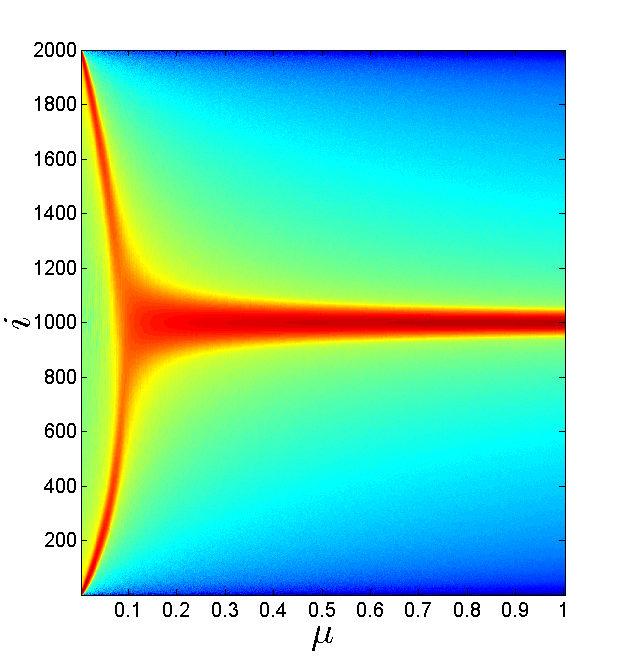}%
  \includegraphics[width=.48\textwidth]{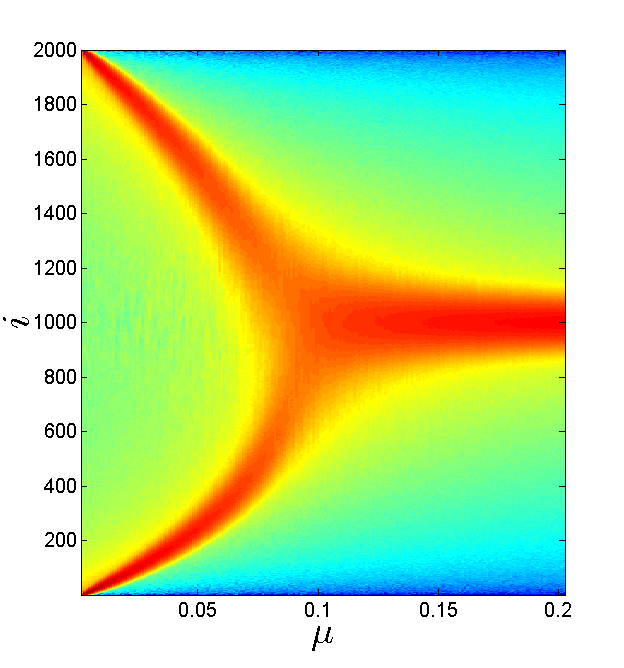}
  \caption{Monte Carlo simulation of the invariant distribution of the
    Markov process with $N=2000$.  Both frames have parameters
    $a=4,b=1,c=3,d=2$; the right frame is  a blowup of the left.}
  \label{fig:4132}
\end{centering}
\end{figure}

\begin{figure}[ht]
\begin{centering}
  \includegraphics[width=.48\textwidth]{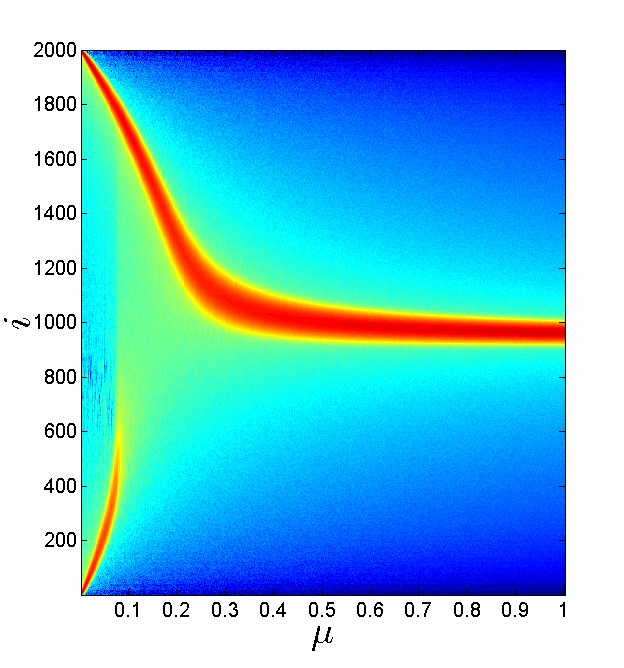}%
  \includegraphics[width=.48\textwidth]{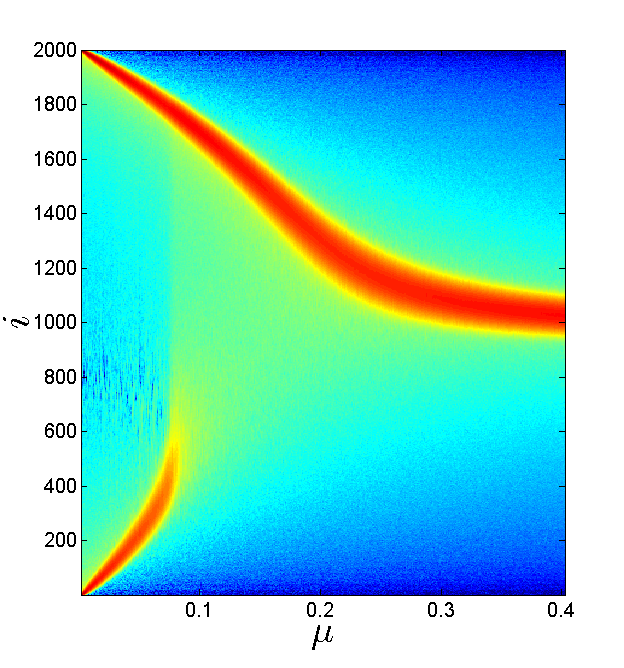}
  \caption{Monte Carlo simulation of the invariant distribution of the
    Markov process with $N=2000$.  Both frames have parameters
    $a=4,b=2,c=1,d=4$; the right frame is  a blowup of the left.}
  \label{fig:4214}
\end{centering}
\end{figure}

In Figures~\ref{fig:4132} and~\ref{fig:4214} we plot the invariant
measure of the the Markov process (in fact, we are plotting the
logarithm of this invariant measure for better
contrast). Figure~\ref{fig:4132} corresponds to an example of case 1.1
and Figure~\ref{fig:4214} is an example of case 2.

The simulations were run for $200,000$ iterations with a burn in of
$20,000$ and averaged over $100$ realizations. The colors towards the
beginning of the rainbow represent densities closer to 1, while the
lower colors represent densities closer to 0.

It is apparent from these pictures that when mutation is rare enough,
there are two metastable mixtures that correspond to the fact that
both strategies are ESS in the absence of mutation.  As the mutation
probability increases, these mixtures approach each other.  In
Figure~\ref{fig:4132}, they seem to merge smoothly into a single
stable mixture, whereas in Figure~\ref{fig:4214} one mixture
disappears suddenly and the other takes over.  We point out that these
mixtures are metastable since we expect switching between them due to
stochastic effects.

\subsection{Connection to Prisoner's Dilemma dynamics}

Above, we have presented the payoff matrix~\eqref{eq:payoff} without
any derivation.  Here we discuss a bit about how such a matrix can
arise.  A common game-theoretic model to consider is that of the
Prisoner's Dilemma (PD), namely a payoff matrix
\begin{equation} \label{eq:pd}
\left[\begin{array}{c|cc}
    & C & D \\\hline
    C & r & s \\
    D & t & p
\end{array}\right],
\end{equation}
where $t > r > p > s$, and we use the same convention of symmetry as
above, i.e. that player2's payoff matrix is the transpose of player1's.
Each player has the option to ``cooperate'' or to ``defect''.  It is
clear that if we assume the strategies are deterministic, the Nash
equilibrium strategy is for both players to defect, since the second
row is componentwise larger than the first.  This is considered a
dilemma because $r>p$, so that if the players could somehow agree to
cooperate and trust each other, they would each receive a higher
payoff, but in any situation where there is a lack of trust they are
better off both defecting.  

Now, we first note that the Prisoner's Dilemma is not in the form
mentioned above, since for us to have two strategies that are both
ESS, we would need that $r>t$ to make $C$ stable.  However, let us
assume that the payoff matrix~\eqref{eq:payoff} is the aggregate
payoff of an $m$-round game between two strategies that cooperate and
defect in some pattern, is it possible to have both strategies be ESS?

\newcommand{\sign}{{\mathrm{sign}}}

The simplest case is if we assume that each strategy is to play $C$
and $D$ with some fixed probabilities, independent of the history of
one's opponent.  For example, if we assume that strategy $A$ is to
cooperate with probability $\alpha\in[0,1]$ and defect with
probability $(1-\alpha)$, and similarly for strategy $B$ with
probability $\beta$, then after a large algebraic computation we find
that $\sign(a-c) = \sign(\beta-\alpha)$ and $\sign(d-b) =
\sign(\alpha-\beta)$, and clearly it is not possible that both of
these be positive.  In particular, a strategy is ESS iff it is more
likely to defect than the opposing strategy.

However, it is possible to construct strategies based on aggregate PD
that allow for memory effects that give multiple ESS.  For example,
consider two strategies known as ``tit-for-tat'' (TFT) and ``always
defect'' (AllD).  In the former, the strategy is to cooperate in the
first round, then do whatever the opponent did in the previous round;
the second strategy is to always defect.  If we assume that strategy
$A$ is TFT and $B$ is AllD, then after $m$ rounds of play, the
aggregate payoffs are
\begin{equation} \label{eq:tftalld}
\left[\begin{array}{c|cc}
      & A & B \\\hline
    A & mr & s+(m-1)p \\
    B & t+(m-1)p & mp
\end{array}\right].
\end{equation}
To see this, note that two players of type $A$ will always cooperate
and two of type $B$ will always defect.  If TFT plays AllD, then the
first player will cooperate once and then defect forever, so will
receive one payoff of $s$ and $m-1$ payoffs of $p$, while the second
player will receive one payoff of $t$ and $m-1$ payoffs of $p$.  Since
$p>s$, we see that $B$ is ESS, but for $A$ to be ESS, we need $mr >
t+(m-1)p$.  But notice that $r>p$, so clearly for $m$ sufficiently
large, $A$ is also ESS.

The types of mixtures that give rise to multiple ESS in a
two-population game have a long history and have been studied
extensively, see~\cite{may1975nonlinear, Nowak.book,
  Imhof.Fudenberg.Nowak.07, wang2008phase, arnoldt2012frequency}.


\section{Deterministic Bifurcation Analysis} \label{sec:deterministic}

The expected step size of the Markov process is given by
\begin{equation*}
  \frac{d\E[X(t)]}{dt} = \omega_+(X_t; \mu) - \omega_-(X_t; \mu),
\end{equation*}
where $\omega_\pm$ is defined
in~(\ref{eq:transitionRate+},\ref{eq:transitionRate-}).  Switching to
the scaled variables and rates 
\begin{equation*}
  x(t) = X_t/N \quad  \omega_\pm(X_t;\mu) = N\Omega_\pm(X_t/N;\mu),
\end{equation*}
we see that the deterministic equation
\begin{equation} \label{eq:deterministic}
\frac{dx}{dt} = f(x;\mu) := \Omega_+(x;\mu) - \Omega_-(x;\mu)
\end{equation}
describes the mean of the Markov process.  As it turns out, in the
limit $N\to\infty$,~\eqref{eq:deterministic} also gives a good
approximations to all paths of the Markov process in a sense that can
be made precise~\cite{Darling.Norris.08}; we discuss this further in
Section~\ref{sec:stochastic}, but for now we
take~\eqref{eq:deterministic} as the system to consider.

The equation~\eqref{eq:deterministic} is a 1-parameter family of
differential equations. We seek to use bifurcation analysis to
determine the equilibria of the deterministic equation. Each case
outlined in section~\ref{ssec:cases} leads to a different class of
bifurcation diagrams (using the mutation rate $\mu$ as the bifurcation
parameter).  The assumption of multiple ESS made above, means that for
$\mu=0$, we have that $x=0, 1$ are nondegenerate fixed points.  Thus
it follows that the system has multiple stable equilibria for some
open set of $\mu$ containing zero.

\subsection{Case 1.1 --- $a+b = c+d$ with $a>d$ and $b<c$} 

Making a substitution for one of the variables allows
$f(x;\mu)$ to be factored into
\begin{align*}
&f(x;\mu) = \bigg(x-\frac{d-b}{a-b-c+d}\bigg)\times\\
&\times\bigg(\frac{-(a-b-c+d)x^2 + (a-b-c+d - (a-b+c-d)\mu)x - 2d\mu}{(a-b-c+d)x^2 + (b+c-2d)x+d} \bigg)
\end{align*}

It is then easy to see that equilibria of the deterministic equation \eqref{eq:deterministic} occur at the following points,
\begin{align*}
x_-(\mu) &= \frac{1}{2} + \bigg(\frac{d-a}{2(d-b)}\bigg)\mu - \frac{\sqrt{4(b-c)^2\mu^2+8(b-d)(a+d)\mu + 4(b-d)^2}}{4(d-b)} \\
x_0(\mu) &= \frac{d-b}{a-b-c+d} = \frac{1}{2} \\
x_+(\mu) &= \frac{1}{2} + \bigg(\frac{d-a}{2(d-b)}\bigg)\mu + \frac{\sqrt{4(b-c)^2\mu^2+8(b-d)(a+d)\mu + 4(b-d)^2}}{4(d-b)}
\end{align*}

\begin{figure}[ht]
\begin{centering}
  \includegraphics[width=.48\textwidth]{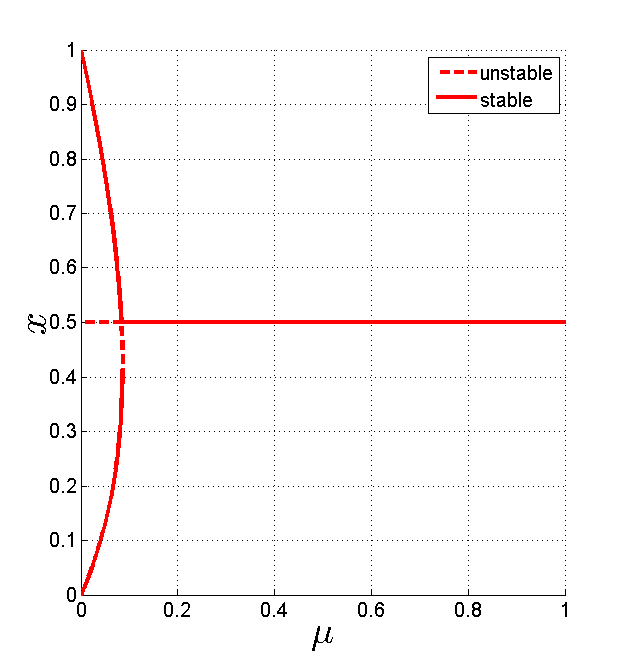}%
  \includegraphics[width=.48\textwidth]{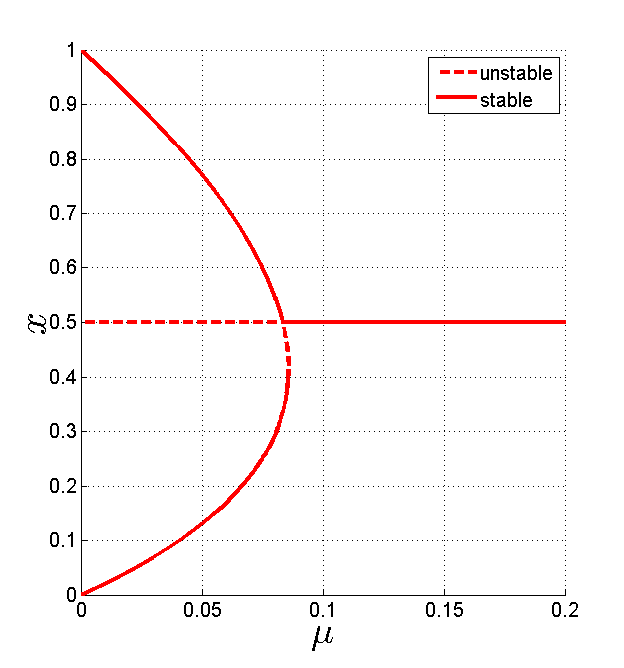}
  \caption{Bifurcation diagram for case 1.1, $a=4,b=1,c=3,d=2$.  Left
    frame is $\mu\in(0,1)$ and the right frame is a blowup.  Compare
    this figure to Figure~\ref{fig:4132}.}
  \label{fig:BD4132}
\end{centering}
\end{figure}

The stability of these equilibria depend on the value of $\mu$. When
$\mu=0$ the equilibria are $x_-(0) = 0$, $x_0(0) = 1/2$, and $x_+(0) =
1$ with $x_-(0)$ and $x_+(0)$ being stable and $x_0(0)$ unstable. The
equilibria maintain these stabilities until the first bifurcation
point, $\mu_1$. This bifurcation point occurs when either $x_-(\mu)$
or $x_+(\mu)$ intersects $x_0(\mu)$. Assuming $a>d$ and $b<c$, this
intersection occurs between $x_+(\mu)$ and $x_0(\mu)$ resulting in
$x_+(\mu)$ becoming unstable and $x_0(\mu)$ becoming stable. Solving
for this intersection, we find $\mu_1 = \frac{d-b}{2(a+d)}$. The
second bifurcation value, $\mu_2$, occurs when $x_-(\mu)$ and
$x_+(\mu)$ become complex-valued leaving $x_0(\mu)$ as the only stable
equilibrium. We find $\mu_2 =(d-b)(a+d-2\sqrt{ad})/(d-a)^2$. (Note
that in case 1.3, there is complete symmetry and $\mu_1 = \mu_2$).
This analysis completely determines the bifurcation diagrams for case
1.1. The bifurcation diagram will have the same qualitative behavior
as Figure \ref{fig:BD4132}.

\subsection{Case 2 --- $a+b>c+d$}

It is no longer possible to factor $f(x;\mu)$, but we can use
the Implicit Function Theorem to describe the motion of the
equilibria.  When $\mu = 0$, equilibria occur at
\begin{equation*}
x_-(0) = 0 ,\quad 
x_0(0) = \frac{d-b}{a-b-c+d} ,\quad 
x_+(0) = 1,
\end{equation*}
again with with $x_-(0)$ and $x_+(0)$ being stable and $x_0(0)$
unstable.  If $x(\mu)$ is an equilibrium of the deterministic equation
then the Implicit Function Theorem implies

$$\frac{dx(\mu)}{d\mu} = - \frac{\dfrac{\partial f}{\partial \mu} (x; \mu)}{\dfrac{\partial f}{\partial x}(x; \mu)}$$

The only
real root of $\dfrac{\partial f}{\partial \mu} (x; \mu)$
between $0$ and $1$ is $$x_* = \frac{-(b-c+2d)+
  \sqrt{(b-c)^2+4ad}}{2(a-b+c-d)}.$$
By observing that
\begin{align*}
\frac{\partial f}{\partial \mu}(x_-(0);0) = 1 >0 ,\quad 
\frac{\partial f}{\partial \mu}(x_+(0); 0) = -1 < 0 \\
\end{align*}
we find that $\dfrac{\partial f}{\partial \mu}(x)>0$ for $x <
x_*$ and $\dfrac{\partial f}{\partial \mu}(x) < 0$ for $x>
x_*$.  The roots of $\dfrac{\partial f}{\partial x}(x, \mu)$
occur when 
\begin{equation*}
  \mu_*(x) = \frac{p(x)}{q(x)},
\end{equation*}
where
\begin{align*}
\end{align*}

By observing that
\begin{align*}
\frac{\partial f}{\partial x}(0, 0) &= \frac{b}{d}-1  < 0 \\
\frac{\partial f}{\partial x}(x_0(0), 0) &= -\frac{(a-c)(b-d)}{ad-bc} > 0 \\
\frac{\partial f}{\partial x}(1, 0) &= \frac{c}{a}-1  < 0
\end{align*}

and that $\mu_*(x)$ is concave down we find that $\mu_*$ is a parabola
in the $\mu$-$x$ plane. Therefore $\frac{\partial f}{\partial
  x}(x, u) >0$ for $\mu < \mu_*(x)$ and $\frac{\partial
  f}{\partial x}(x, u) < 0$ for $\mu > \mu_*(x)$.

\begin{figure}[th]
\begin{center}
\includegraphics[scale=0.4]{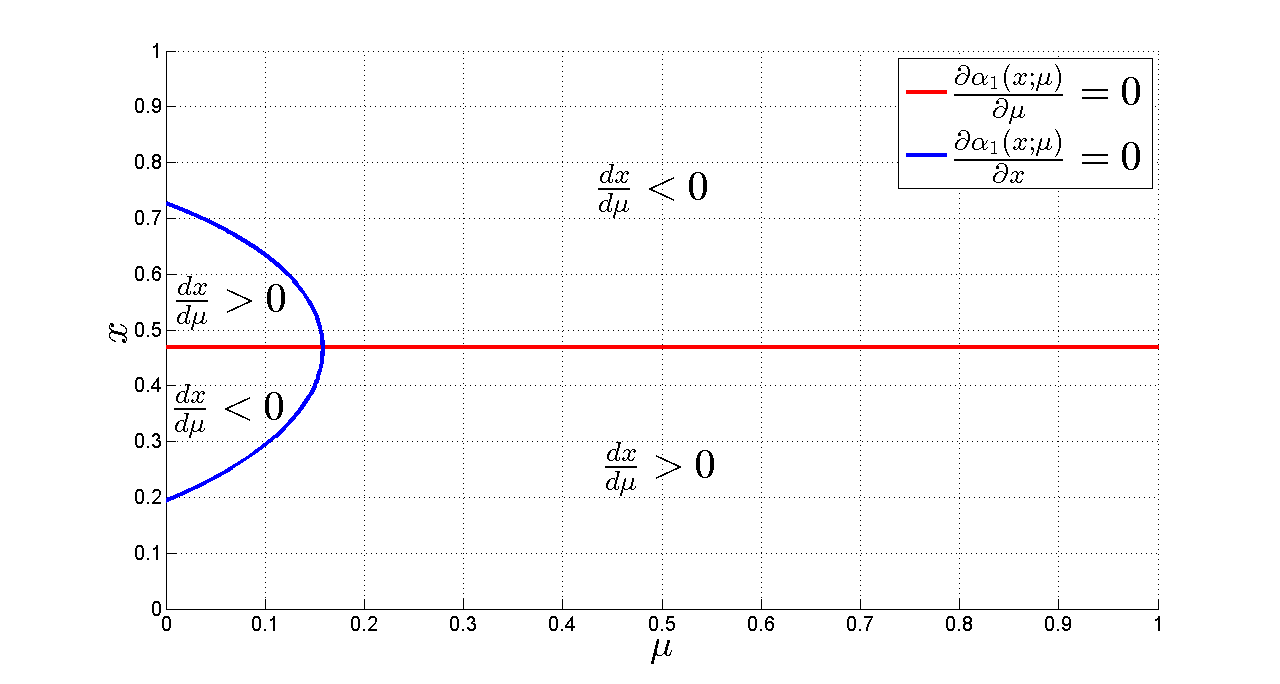}
\caption{Sign of $dx/du$ for $a+b > c+d$}
\label{fig:dxdu}
\end{center}
\end{figure}

\begin{figure}[ht]
\begin{centering}
  \includegraphics[width=0.49\textwidth]{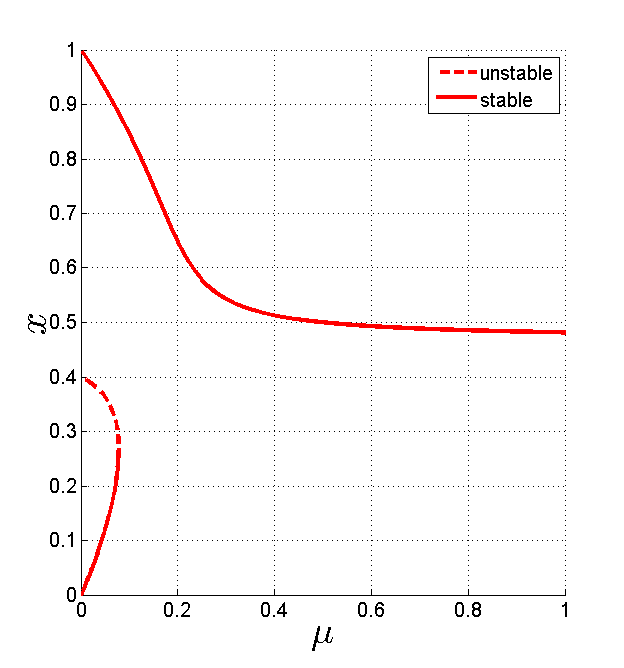}%
  \includegraphics[width=0.49\textwidth]{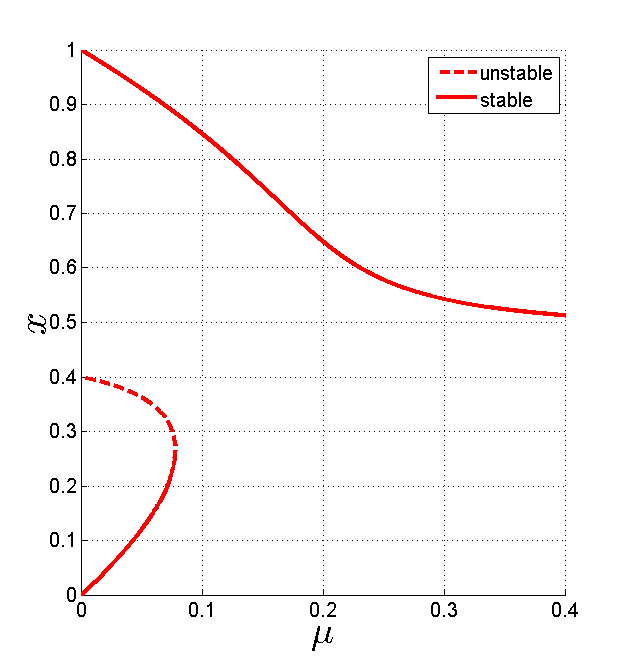}
  \caption{Bifurcation Diagram of case 2: $a=4, b=2, c=1, d=4$}
  \label{fig:BD4214}
\end{centering}
\end{figure}

Putting the information about $\frac{\partial f}{\partial \mu}(x)$ and
$\frac{\partial f}{\partial x}(x, \mu)$ together, there are $4$
regions where $\frac{dx}{d\mu}$ changes sign. The first equilibrium
$x_-(\mu)$ starts at $0$ which is located in a region where
$\frac{dx}{d\mu} > 0$ and is therefore increasing. It is easy to show
that $x_0(0) < x_*$ and thus $x_0(\mu)$ starts in a region where
$\frac{dx}{d\mu} < 0$. Therefore $x_0(\mu)$ is always
decreasing. Similarly $x_+(\mu)$ starts in a region where
$\frac{dx}{d\mu} < 0$ and therefore is always decreasing. A
bifurcation point occurs at the value of $\mu$ when two of the
equilibria become complex valued. By continuity, this must occur along
$u_*(x)$ when $x_+(\mu)$ and $x_0(\mu)$ collide. The bifurcation
diagram described by the motion of these equilibria has the
qualitative behavior as Figure~\ref{fig:BD4214}.

\section{Stochastic Analysis} \label{sec:stochastic} Section
\ref{sec:deterministic} described the deterministic behavior of the
process in the $N\to\infty$ limit.  However, if $N\gg 1$ but finite,
the Markov chain model with retain stochasticity, although of course
 the variance of this process will vanish as $N\to\infty$.

In this section, we study this random process for $N$ large but
finite, using the standard tools~\cite{Gardiner.book}.  The two
main asymptotic expansions used here are the ``diffusion
approximation'' and the ``WKB expansion'', which roughly go as
follows.  In both cases, we assume that we start off with a Master
equation for the stochastic process after a system size scaling,
specifically we write
\begin{equation} \label{eq:masteri}
\frac{\partial P(i, t)}{\partial t} = \omega_+(i-1)P(i-1, t) + \omega_-(i+1)P(i+1, t) - (\omega_+(i) + \omega_-(i))P(n, t)
\end{equation}
with reflecting boundary conditions at $i=0$ and $i=N$.  In terms of
the scaled variable $x=i/N$, this master equation is
\begin{equation}\label{eq:masterx}
\begin{split}
  \frac{\partial P(x, t)}{\partial t} 
  &= N \Omega_+(x-1/N) P(x-1/N, t) + N \Omega_-(x+1/N) P(x+1/N, t) \\
  &\quad- N(\Omega_+(x) + \Omega_-(x)) P(x, t).
\end{split}
\end{equation}

In the diffusion approximation, one simply makes a Taylor series
expansion of~\eqref{eq:masterx} and then truncates at second order,
which gives a Fokker--Planck PDE.  One can then consider the paths of
this Fokker--Planck equation which are governed by an SDE.  Computing
the stationary distribution and the jumping timescale then corresponds
to solving this PDE with particular boundary conditions.  In
particular, since the system is one-dimensional, it can be written in
the form of a diffusion in a one-dimensional potential $\Phi$.  This
potential governs the dynamics of this system, and in particular
governs the switching times of the system from one basin to another.

In the WKB approximation, the procedure is more complex.  We seek to
understand how the distribution of an initial condition placed near
one attractor will behave.  The motivating presumption is that it will
settle down to a characteristic shape inside the basin of attraction
of this point, and this shape will decay to zero on a timescale
associated with the escape of trajectories out of this basin.  More
concretely: We first assume that the solution of the Master equation
assuming that the initial condition is at $x_-$ w.p.1 can be written
in {\em quasistationary form},
\begin{equation*}
  \Pi(x,t)  =e^{-t/\tau_-} \exp\left(N \Psi_{-1}(x) + \Psi_0(x) + N^{-1}\Psi_1(x)+\cdots\right),
\end{equation*}
where $\tau_-$ is the yet-to-be-determined jumping timescale.  If we
truncate the above expression to remove terms of order $N^{-1}$, we
obtain the quasistationary distribution
\begin{equation*}
  \widehat\Pi(x) := e^{t/\tau_-} \Pi(x,t) = k(x) e^{-N\Psi_{-1}(x)},\quad \log k(x) = \Psi_0(x).
\end{equation*}
We drop the subscript on $\Psi_{-1}$ in the sequel.  Here the function
$\Psi$ plays the role of a potential, similar to $\Phi$ in the
diffusion approximation and is also called the {\em quasipotential}
for the system.  The distribution $\widehat\Pi(x)$ is not a stationary
distribution, since it will decay on the long timescale $\tau_-^{-1}$,
but the idea is that the solution relaxes very quickly to
$\widehat\Pi(x)$, and then decays slowly to zero.  One can perform an
asymptotic expansion in inverse powers of $N$ to obtain equations for
$k,\Psi, \tau_-$, see details below.

It is known~\cite{Bressloff.book, Bressloff.10} that the WKB expansion
gives an asymptotically correct expression for the jumping timescale,
whereas the diffusion approximation gives a slightly different, and
wrong expression, but we observe two things for this problem: first,
the diffusion approximation is much simpler to compute, and second,
the answers obtained from the numerics show that the two answers are
indistinguishable on the types of problems that we can simulate.

\subsection{Diffusion Approximation} \label{sec:diffusion}

Performing a Taylor expansion on~\eqref{eq:masterx} in powers of $1/N$
results in
\begin{equation} \label{eq:masterTaylor}
\frac{\partial P(x, t)}{\partial t} = \sum_{n=1}^{\infty} N^{1-n} \frac{(-1)^n}{n!} \frac{\partial^n}{\partial x^n} \left[ \alpha_n(x) P(x) \right]
\end{equation}
where $\alpha_n(x) = \Omega_+(x) + (-1)^n \Omega_-(x)$.  Terminating
this series to order $1/N$ and writing $f= \alpha_1, \sigma =
\alpha_2$ results in the {\bf diffusion approximation}
\begin{equation} \label{eq:diffusion}
\frac{\partial P(x, t)}{\partial t} = -\frac{\partial}{\partial x} \left[ f(x) P(x) \right] + \frac{1}{2N} \frac{\partial^2}{\partial x^2} \left[ \sigma(x) P(x) \right] =: \mcL^* P
\end{equation}
The diffusion approximation is a Fokker-Planck equation and has the equivalent SDE form
\begin{equation}\label{eq:SDE}
dx = f(x) \ dt + \sqrt{\frac{\sigma(x)}{N}} \ dW(t).
\end{equation}
As $N \rightarrow \infty$, this is exactly the deterministic equation
\eqref{eq:deterministic} discussed in section \ref{sec:deterministic}.


\subsubsection{Linear noise Approximation} \label{sssec:vk}

The linear noise approximation, more commonly known as the van Kampen
expansion~\cite{vanKampen.book, vanKampen.82, Gardiner.book}.  The idea here is to
write $x(t) = \varphi(t) + N^{-1/2}z$, obtain equations for $z$ , then
assume that $\varphi(t)$ has some well-defined limit as $t\to\infty$.
In the case where the deterministic system has a unique attracting
fixed point $x^*$, applying the van Kampen expansion (see
Appendix~\ref{sec:vK} for more detail, we obtain
\begin{equation}\label{eq:VKmoments}
  \E[x(t)] \to x^*,\quad \Var(x(t)) \to \frac{ \sigma(x_*)}{2 N \av{f'(x_*)}}.
\end{equation}
This computation is simple and gives useful information: in particular
we obtain the scaling of the variance as a function of $N$.  The one
issue with this approximation is that it is guaranteed to be Gaussian
and thus cannot capture the skew of the original process.

\subsubsection{Higher-order Moments} \label{sec:perturbation} The
linear noise approximation is the lowest order expansion for the
diffusion approximation, but this can be extended to a general
perturbation expansion on~\eqref{eq:SDE}.  Writing
\begin{equation*}
  x(t) = x_* + N^{-1} x_1(t) + N^{-2}x_2(t)+\dots
\end{equation*}
and truncating at second order, we obtain for large $t$ (see
Section~\ref{app:ho} for details):
\begin{align}
&\E[x(t)] \to x_* + \frac{f''(x_*)\sigma(x_*)}{4N(f'(x_*))^2}, \label{eq:moment1} \\
&\Var(x(t)) \to \frac{\sigma(x_*)}{2N\av{f'(x_*)}} + \frac{(f''(x_*))^2 \sigma^2(x_*)}{8N^2(f'(x_*))^4} + \frac{(\sigma'(x_*))^2}{16N^2(f'(x_*))^2}, \label{eq:moment2} \\
&\E[(x(t) - \E[x(t)])^3] \to \frac{(f''(x_*))^3 \sigma(x_*)^3}{8N^3f'(x_*)^6} + \frac{3(\sigma'(x_*))^2 f''(x_*) \sigma(x_*)}{32 N^3 f'(x_*)^4}. \label{eq:moment3}
\end{align}
This next approximation gives additional insight into our process for
case 1. In section \ref{sec:deterministic}, we determined that once
$\mu \geq \mu_2$, there is one stable equilibria at $x_0 = 1/2$. At
this equilibrium, both strategies have the same proportions, but the
sign of the centered third moment can give insight into which strategy
has a slight advantage.  For now, we will just determine the sign of
the centered third moment. Around the equilibria $x_0 = 1/2$, equation
\eqref{eq:moment3} simplifies to
\[\E[(x - \E[x])^3] = \frac{(f''(x_0))^3}{64N^3f'(x_0)^6}\]
The sign of this only depends on the sign of $f''(x_0)$. Using
the fact that $a+b=c+d$ we can factor $f''(x_0)$ as 
\[ f''(x_0) = \frac{32(c+d)(a-c)(a-d)(2\mu - 1)}{(a + b + c + d)^3} \]

The sign of the skew is now easily found.

\[ \begin{array}{c|cc}
 & \mu_2 <\mu < .5 & \mu>.5 \\ \hline
a>d & \E[(x - \E[x])^3] < 0 & \E[(x - \E[x])^3] > 0 \\
a<d & \E[(x - \E[x])^3] > 0 & \E[(x - \E[x])^3] < 0
\end{array} \]

\subsubsection{Stationary Distribution}
Setting the left-hand side of equation \eqref{eq:diffusion} to 0 and
using the reflecting boundary conditions, it is straightforward to
find the stationary distribution for the diffusion approximation:
\begin{equation}\label{eq:stationary3}
P_s(x) = \frac{\mathcal{B}e^{-N\Phi(x)}}{ \sigma(x)}
\end{equation}
where $\mathcal{B}$ is a normalization constant and the {\bf
  quasipotential} if the function $\Phi(x) = -2\int_0^x \frac{
  f(x')}{ \sigma(x')} \ dx'$.  By differentiating
equation~\eqref{eq:stationary3} and neglecting the $1/N$ terms, one
sees that for large $N$, the maxima of the stationary distribution
occur at the stable fixed points of the deterministic equation
\eqref{eq:deterministic}, while the minima of the stationary
distribution occur at the unstable fixed points of the deterministic
equation~\eqref{eq:deterministic}.


\subsubsection{First-Passage Times using Diffusion}

Here we see to compute the mean first passage time (MFPT) from $x_-$
to $x_+$ using~\eqref{eq:diffusion}.

To calculate the MFPT from $x_-$ to $x_+$ we restrict our attention to
the interval $[0, x_+]$ and impose a reflecting boundary conditions at
$0$ and absorbing boundary conditions at $x_+$. Let $T(x)$ denote the
first passage time (FPT) at which the system reached $x_+$ starting at
$x$ and $\tau_-(x) = \E[T(x)]$ denote the MFPT from $x$ to
$x_+$.  The MFPT satisfies
\begin{equation}\label{eq:odetau}
  \mcL \tau_- +1=0, \quad \tau_-'(0) = 0, \quad  \tau_-(x_+) = 0,
\end{equation}
where
\begin{equation*}
  \mcL  := f(x) \frac{\partial
  }{\partial x} + \frac{\sigma(x)}{2N} \frac{\partial^2
  }{\partial x^2}
\end{equation*}
is the (formal) adjoint of the Fokker--Planck operator defined
in~\eqref{eq:diffusion}.  This follows from the fact that $\mcL$ is
the generator of the diffusion process~\cite{Gardiner.book}.  Since
the coefficients of the diffusion are known explicitly, we can obtain
a formula for $\tau_-(x)$ by quadrature and apply Laplace's
integral method (see Appendix~\ref{app:Laplace} for details) to obtain
\begin{equation} \label{eq:tau-}
\tau_- = \frac{2 \pi}{\Omega_+(x_-) \sqrt{\Phi''(x_-) | \Phi''(x_0)|}} e^{N(\Phi(x_0)-\Phi(x_-))}.
\end{equation}
Notice in particular that this approximation is independent of $x$ as
long as it is not too close to $x_0$.  Thus we can think of the
quantity $\tau_-$ as the MFPT for any initial condition
starting in a neighborhood of $x_-$ to reach a neighborhood of $x_+$.

We can perform the same computation in the other direction to obtain
\begin{equation} \label{eq:tau+}
\tau_+ = \frac{2 \pi}{\Omega_-(x_+) \sqrt{\Phi''(x_+) | \Phi''(x_0)|}} e^{N(\Phi(x_0)-\Phi(x_+))}.
\end{equation}


\subsection{ WKB Approximation} \label{sec:WKB}

\newcommand{\Pirel}{\Pi_{\mathsf{rel}}}
\newcommand{\Piact}{\Pi_{\mathsf{act}}}

We seek a quasistationary distribution of the master equation
\eqref{eq:masterx} in the form
\begin{equation} \label{eq:defofPi}
\Pi(x,t) = k(x) e^{-N\Psi(x)} e^{-t/\tau_-}.
\end{equation}
If we plug this Ansatz into the Master equation~\eqref{eq:masterx},
Taylor expand $k,\Psi$ in powers of $N^{-1}$ (see
appendix~\ref{app:WKB} for details), and assume that $\tau_- \gg N$,
then we obtain two solutions for the WKB solution, called the
relaxation and activation solutions:
\begin{align} 
\Pirel(x) &= \frac{B}{\Omega_+(x) - \Omega_-(x)},\label{eq:relaxation}\\
\Piact(x) &= \frac{A}{\sqrt{\Omega_+(x)\Omega_-(x)}} e^{-N\Psi(x)},\quad \Psi(x) = \int_0^x \ln\left( \frac{\Omega_-(x')}{\Omega_+(x')}\right) \ dx'..\label{eq:activation}
\end{align}
These solutions have a total of two free constants to be determined.
Performing a matched asymptotic expansion around the point $x_0$
(again see appendix~\ref{app:WKB} for details), we obtain the formula
\begin{equation}\label{eq:MFPT1}
\tau_- = \frac{2\pi}{\Omega_+(x_-)\sqrt{|\Psi''(x_0)|\Psi''(x_-)}} e^{N(\Psi(x_0)-\Psi(x_-))}
\end{equation}
Similarly, the MFPT from $x_+$ to $x_-$ is
\begin{equation}\label{eq:MFPT2}
\tau_+ = \frac{2\pi}{\Omega_-(x_+)\sqrt{|\Psi''(x_0)|\Psi''(x_+)}} e^{N(\Psi(x_0)-\Psi(x_+))}
\end{equation}


\subsection{Comparison of quasipotentials}\label{sec:compare}

It is useful to compare the expressions~\eqref{eq:tau-}
and~\eqref{eq:MFPT1}; notice that the only difference is that there is
a different quasipotential: $\Phi$ for the diffusion approximation and
$\Psi$ for the WKB approximation.  We repeat these here for comparison:
\begin{equation*}
  \Phi(x) = \int_0^x -2\frac{\Omega_+(x') - \Omega_-(x')}{\Omega_+(x') + \Omega_-(x')}\,dx',\quad \Psi(x) = \int_0^x \log\left(\frac{\Omega_-(x')}{\Omega_+(x')}\right)\,dx'.
\end{equation*}
Now first notice that the only way in which $\Phi,\Psi$ appear in the
timescale is through the difference of values at $x_0$ and $x_-$, so
in fact we are really concerned with the numbers we obtain when we
evaluate these integrals from $x_-$ to $x_0$.

These integrands are certainly different functions and there is no
reason why they should be close.  However, let us write $q(x) =
\Omega_-(x)/\Omega_+(x)$, and expand each integrand in a Taylor series
at $q=1$:
\begin{equation*}
\Phi'(x) = -2\frac{1-q}{1+q}= \sum_{k=1}^\infty \frac{(-1)^{k-1}}{2^{k-1}} (q-1)^k,\quad   \Psi'(x) = \log(q) = \sum_{k=1}^\infty \frac{(-1)^{k-1}}{k} (q-1)^k
\end{equation*}

Note that these expressions are equal up to the third order term in
$(q-1)$, and thus as long as $\Omega_+(x)$ and $\Omega_-(x)$ are
close, then we expect these two integrals to be close.

\begin{figure}[ht]
\begin{centering}
  \includegraphics[width=.45\textwidth]{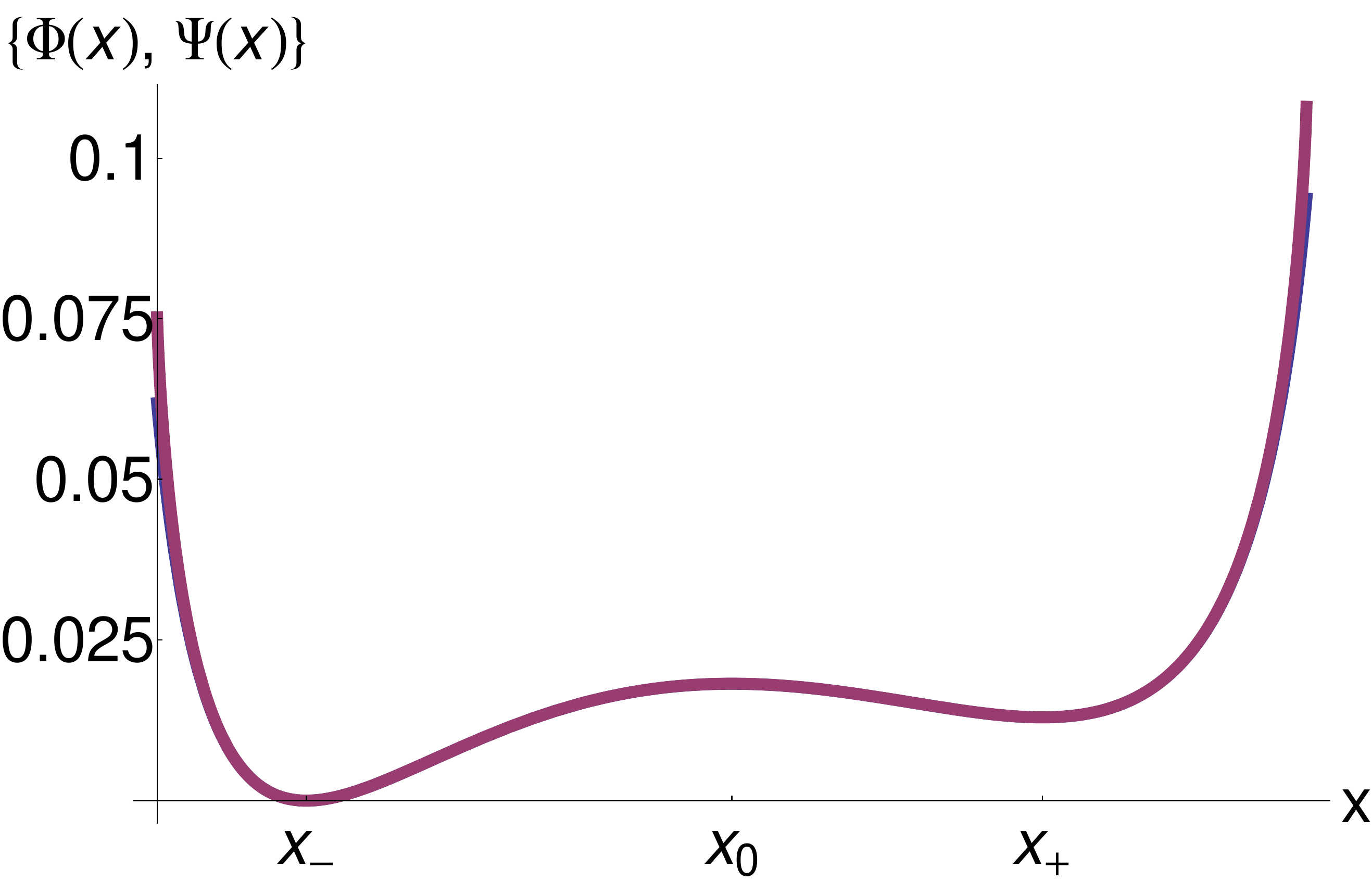}%
  \includegraphics[width=.45\textwidth]{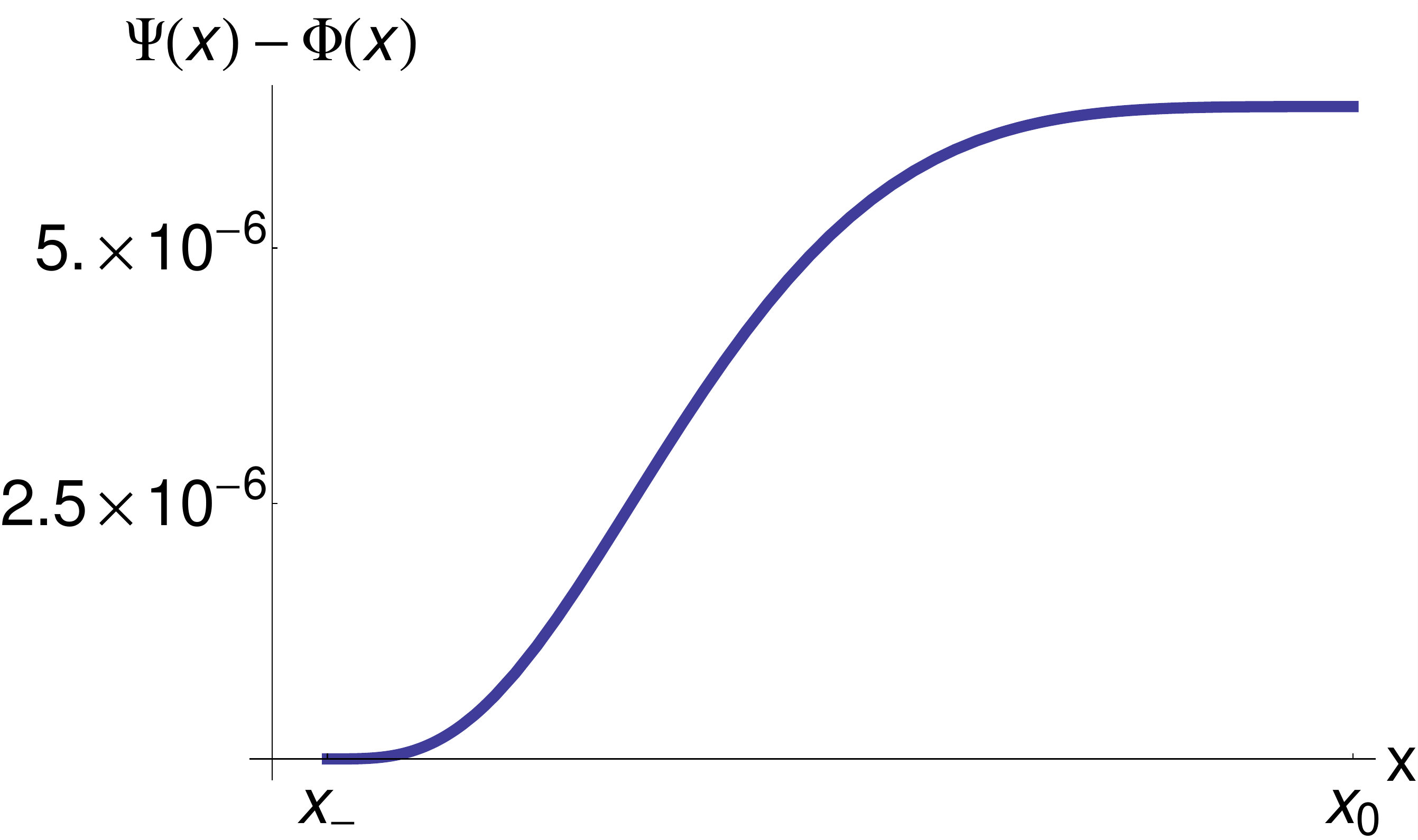}
  \caption{Comparison of the quasipotentials $\Phi(x),\Psi(x)$ for
    particular parameter values $a=4,b=1,c=3,d=2,\mu=0.05$.  The left
    frame shows the two functions on the same axis, and the right
    frame shows their difference.}
  \label{fig:quasi}
\end{centering}
\end{figure}

In particular, we plot the two quasipotentials for the parameters
$a=4,b=1,c=3,d=2,\mu=0.05$ in Figure~\ref{fig:quasi}.  Note that the
two functions are indistinguishable except for near the ends of the
domain.  Considering the difference on the interval $[x_-,x_0]$, and
we see that the difference between the functions is four orders of
magnitude smaller than the functions themselves.

Now, of course, when $N\to\infty$, the quasipotentials are multiplied
by $N$ and exponentiated, so that even a small difference in the
quasipotential will make a huge difference asymptotically.  However,
we should note that both approximations are only expected to be
accurate in the limit $N\to\infty$, and for any finite but large $N$,
the neglected terms in the expansion could dominate any difference in
the quasipotentials.  We actually see that this is so in the
simulations presented below.


\subsection{Simulation Results}
%
%
%

\subsubsection{Comparison between van Kampen approximation and
  deterministic system} \label{ssec:comparisonMV} This section
compares the mean and variance of the simulations of the Markov
process to the mean and variance calculated from the linear-noise
approximation in section \ref{sssec:vk}.

\begin{figure}[th]
\begin{center}
\includegraphics[width=0.49\textwidth]{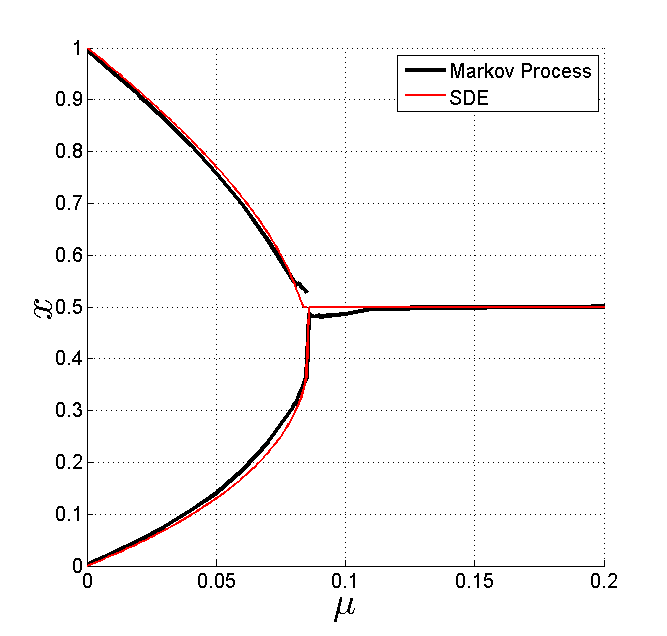}%
\includegraphics[width=0.49\textwidth]{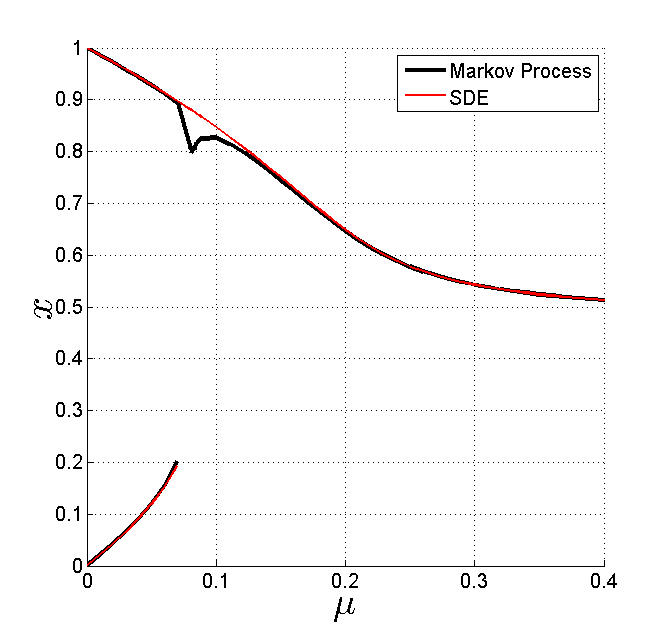}
\caption{Mean of Markov process versus predicted mean.  In the left
  frame we have an example of Case 1.1, $a=4, b=1, c=3, d=2, N=2000$;
  in the right we have $a=4, b=1, c=2, d=4, N=2000$.}
\label{fig:4132mean}
\end{center}
\end{figure}

Recall that the mean of the linear noise approximation is given by the
deterministic
equation~\eqref{eq:deterministic}. Figure~\ref{fig:4132mean} compares
the mean of the simulation of the Markov process against the
bifurcation diagrams determined in section~\ref{sec:deterministic}. We
see that the deterministic equation is a good approximation of the
mean of the process. The only location of disagreement is for $\mu$
close to the bifurcation values. Of course, as $N$ increases this will
become more accurate.

\begin{figure}[th]
\centering
\subfigure[Conditional variance around $x_-$]{%
\includegraphics[scale=.4]{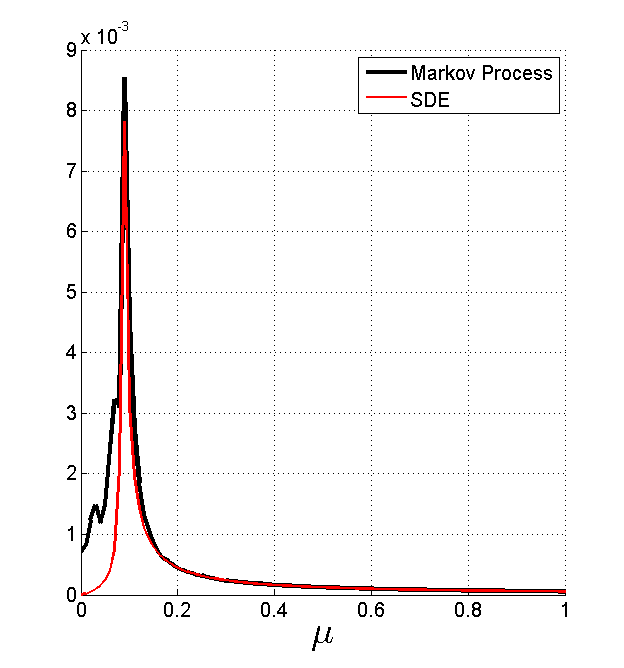}
\label{fig:4132var0}}
\quad
\subfigure[Conditional variance around $x_+$]{%
\includegraphics[scale=.4]{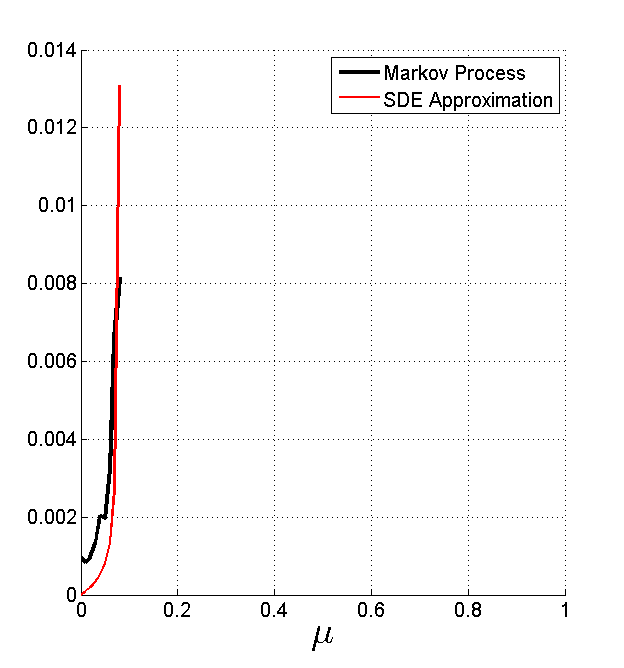}
\label{fig:4132var1}}
\caption{Conditional variance of case 1.1: $a=4, b=1, c=3, d=2,
  N=2000$.}
\label{fig:4132var}
\end{figure}

Figure~\ref{fig:4132var} compares the conditional variance of the Markov
process around each equilibrium point against equation~\eqref{eq:var}.
For the purposes of the figures, when there is only one equilibrium,
we plot it on the $x_-$ figure.  Similar to the mean, the figures
differ slightly once $u$ gets close to the bifurcation points. Again
this is due the size of $N$.

\subsubsection{Comparison of MFPT} \label{ssec:comparison}

The MFPTs calculated from both the SDE approximation and WKB
approximation differ in the prefactor and the form of the
quasipotential.  The MFPT calculated from the WKB approximation will
more accurately capture the exponential time scale in which switching
occurs in the limit as $N\to\infty$.  However, in simulating both
formulations we see they are indistinguishable. Simulations can only
be done for relatively small $N$ due to the exponential stiffness of
the problem, and we see that for any accessible simulation, the
relative error is indistinguishable.  

To give us an idea of the accuracy of these approximations, we use a
Monte Carlo simulation of the Markov process: We averaged 1000
realizations of the Markov process to calculate the MFPT from $x_-$ to
$x_+$. We know from both the SDE approximation~\eqref{eq:tau-} and the
WKB approximation~\eqref{eq:MFPT1} that these MFPTs are exponentially
large, which makes them hard to simulate using the Markov process.

\begin{figure}[th]
\centering
\subfigure[]{%
\includegraphics[width=0.45\textwidth]{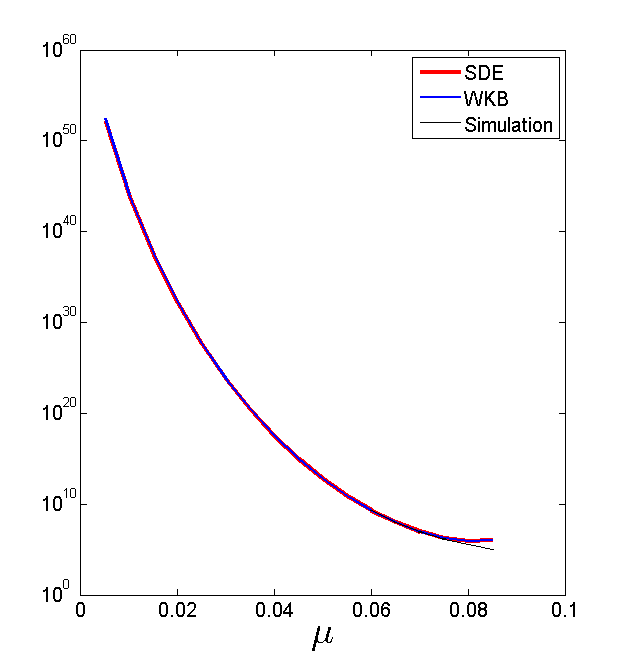}
}
\quad
\subfigure[]{%
\includegraphics[width=0.45\textwidth]{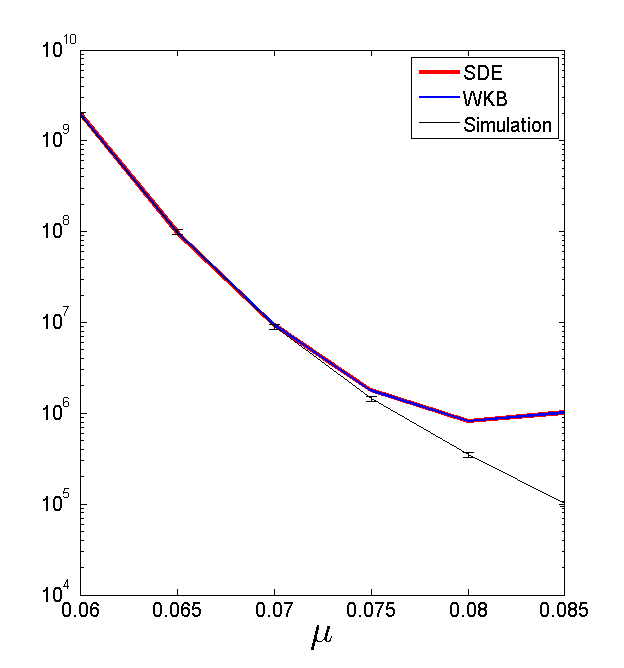}
}
\subfigure[]{%
\includegraphics[width=0.45\textwidth]{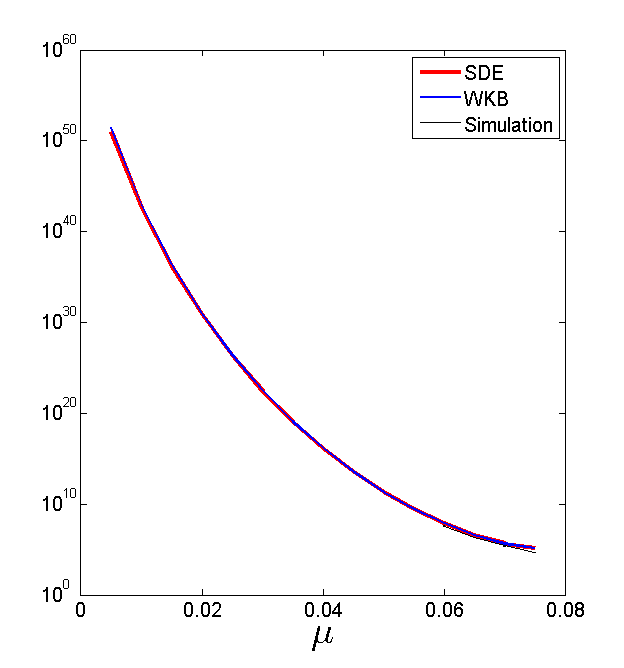}
}
\quad
\subfigure[]{%
\includegraphics[width=0.45\textwidth]{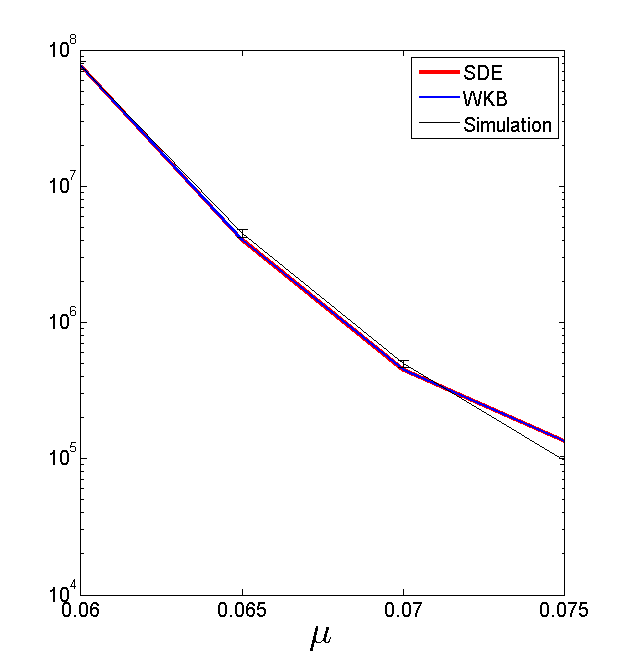}
}
\caption{MFPT from $x_-$ to $x_+$ in two cases: the top row is an
  example of case 1.1: $a=4, b=1, c=3, d=2, N=1000$, whereas the
  bottom row is an example of case 2: $a=4, b=2, c=1, d=4, N=1000$.
  In all cases we are plotting $\tau_-$, the transition time from
  $x_-$ to $x_+$.  The right column is a blowup of the left.}
\label{fig:TM}
\end{figure}

We plot all of these in Figure~\ref{fig:TM}.  As $\mu\to0$, the
transition time becomes exponentially large, so the range of $\mu$ for
which we can capture this transition time is quite restricted.  In
each row, we plot only those $\mu$ where we could perform a Monte
Carlo simulation on the Markov chain directly in the right column.  In
the left column, we extend the plotting region down to small $\mu$.

For case 1.1, we were able to simulate the MFPT starting at $\mu =
.06$. Figure compares the MFPT from $x_-$ to $x_+$ of the simulation
of the Markov process (along with the $95 \%$ confidence interval) to
the MFPT approximated by the SDE and WKB approximations.  Initially
$\tau_-^{SDE}$ and $\tau_-^{WKB}$ are indistinguishable and both fall
within the confidence interval of our simulated $\tau_-$. As $\mu$
approaches the bifurcation point, $\tau_-^{SDE}$ and $\tau_-^{WKB}$
are still indistinguishable, but they no longer fall within the
confidence interval of $\tau_-$. This is due to the fact that
$\tau_-^{SDE}$ and $\tau_-^{WKB}$ are good approximations when $N$ is
large because the stationary distribution is sharply peaked around
$x_-$ and $x_+$. The simulations we ran use $N = 1000$ and as $\mu$
approaches the bifurcation point, the stationary distribution is not
as sharply peaked, and we expect to need larger $N$ to make the
asymptotics accurate.

%
%
%
%
%


\section{Discussion} \label{sec:discussion}

In this paper, we have studied a stochastic version of a standard
game-theoretic model, concentrating on the case where there are
multiple ESS strategies and when the population is small enough that
stochastic effects need to be taken into account.  Our dynamics also
include a mutation parameter for offspring, and we prove that a large
enough mutation factor always causes the multiple ESS solutions to
merge.  The presence of noise due to a finite population turns each
ESS solution into a metastable solution, and we computed the
asymptotics of the switching timescales between these.

We were able to completely characterize the bifurcation diagrams for
this model when the system is in the infinite population limit.  We
also were able to make two nonintuitive observations about the
stochastic process:  

\begin{enumerate}

\item We see that the computation for the switching time using the
diffusion approximation or the asymptotics directly on the Master
equation give essentially the same results: specifically, we see in
Figure~\ref{fig:TM} that the two approximations are indistinguishable
over a wide range of parameters.  What is perhaps more surprising is
that these two agree even when they are far away from the actual
answer gotten from Monte Carlo simulations --- these approximations
are more like each other than the quantity that they are seeking to
approximate.  Of course, the fact that these approximations deviate
from the true value determined by Monte Carlo is not a contradiction,
since these asymptotics are only guaranteed to be correct in the
$N\to\infty$ limit.  We also see from the same figure that this limit
is very non-uniform in system parameters, e.g. for fixed $N$ the
approximation has a different level of effectiveness for different
$\mu$.

\item We also see that the stochastic system gives a finer picture of
  the dynamics than the deterministic limit does.  For example,
  consider cases where $\mu$ is small enough that the system has
  multiple ESS states.  In the deterministic regime ($N=\infty$), we
  have multiple stable solutions that will persist for all time, and
  in particular, the long-term steady state behavior will depend
  strongly on the initial condition of the system.  For large $N$,
  however, the system will spend most of its time near these two
  solutions, but will switch between the two on the timescales
  computed in Section~\ref{sec:stochastic}.  Of course, in general the
  switching timescales in each direction will be different (often,
  significantly different) so this allows us to understand which of
  the two solutions are ``more'' stable, and which solution we will
  expect to observe more often.  We also saw in
  Section~\ref{sec:perturbation} that the stochastic dynamics have a
  bias even when the deterministic dynamics are symmetric, so that the
  finite-size effects can lead to a bias; for example, the stochastic
  dynamics have a skew related to the self-interaction strength of
  each population.
	
\end{enumerate}

\section{Acknowledgments}

L. D. was supported by the National Aeronautics and Space
Administration (NASA) through the NASA Astrobiology Institute under
Cooperative Agreement Number NNA13AA91A issued through the Science
Mission Directorate.  M. G. acknowledges support from National Science
Foundation grant DMS 08-38434 ``EMSW21-MCTP: Research Experience for
Graduate Students.''

\providecommand{\bysame}{\leavevmode\hbox to3em{\hrulefill}\thinspace}
\providecommand{\MR}{\relax\ifhmode\unskip\space\fi MR }
\providecommand{\MRhref}[2]{%
  \href{http://www.ams.org/mathscinet-getitem?mr=#1}{#2}
}
\providecommand{\href}[2]{#2}

\appendix

\section{Details of derivations above}

\subsection{van Kampen approximation}\label{sec:vK}

The linear noise approximation developed by van Kampen~\cite{vanKampen.book,
vanKampen.82, Gardiner.book} assumes the system has linear noise and
further expands the diffusion approximation to obtain the lowest order
expansion. The linear-noise expansion allows us to easily calculate
the moments of this approximation. We begin by assuming
\begin{equation}\label{eq:cov}
x = \phi(t) + N^{-1/2} z
\end{equation}
where $z$ is a fluctuation term and $\phi(t)$ is not yet
determined. Making this change of variables in equation
\eqref{eq:masterTaylor}, requiring $\phi(t)$ to satisfy $\phi'(t) =
f(\phi(t))$ (deterministic equation) and collecting only the
$O(1)$ terms gives
\begin{equation}\label{eq:FPEz}
\frac{\partial P(z, t)}{\partial \tau} =  - f^{\prime}(\phi(t))\frac{\partial}{\partial z} \left[ z  P(z, t) \right] +\frac{1}{2} \sigma(\phi(t)) \frac{\partial^2}{\partial z^2}\left[P(z, t) \right]
\end{equation}
By multiplying equation \eqref{eq:FPEz} by $z$, integrating, and then
integrating by parts the appropriate number of times, we find that $z$
is Gaussian with moment equations,
\begin{equation} \label{eq:firstmomenteq}
\frac{\partial}{\partial t} \E[z] = f'(\phi(t)) \E[z]
\end{equation}
and
\begin{equation} \label{eq:secondmomenteq} \frac{\partial}{\partial t}
  \E[z^2] = 2 f^{\prime}(\phi(t)) \E[z^2] + \sigma(\phi(t)).
\end{equation}

Assuming that if the process runs long enough, the first and second
moments will eventually be stationary, and noting that $\phi(t)\to
x^*$, we find
\begin{equation*}
\E[z] = 0, \quad \E[z^2] = -\frac{\sigma(x_*)}{2 f'(x_*)}.
\end{equation*}
where $x_*$ is an attracting fixed point of the deterministic
equation. Making a change of variable back to $x$ results in,
\begin{equation*}
\E[x] = \phi(t),\quad
\E[x^2] = \phi^2(t) -\frac{\sigma(x_*)}{2 f'(x_*)}.
\end{equation*}

This shows that the mean of the linear noise approximation evolves
according to the deterministic equation and the variance evolves
according to
\begin{equation}\label{eq:var}
\Var(x) = -\frac{ \sigma(x_*)}{2 N f'(x_*)}.
\end{equation}

\subsection{Higher-order corrections}\label{app:ho}
We find higher order corrections to the moments of the linear noise
approximation. We begin by assuming $$x(t) = x_*(t) + \varepsilon
x_1(t) + \varepsilon^2 x_2(t) + \cdots$$ where $\varepsilon =
\frac{1}{N}$.  Plugging the Taylor expansion of $f(x)$ and
$\sqrt{\sigma(x)}$ into the the SDE \eqref{eq:SDE} and equating powers
of $\varepsilon$, we obtain the following system of SDEs which can be
solved iteratively.
\begin{align}
dx_0 &= f(x_*) \ dt, \label{eq:SDE1}\\
dx_1 &= f'(x_*)x_1 \ dt + \sqrt{ \sigma(x_*)} \ dW(t), \label{eq:SDE2}\\
dx_2 &= \left( f'(x_*) x_2 + \frac{f''(x_*) x_1^2}{2} \right) \ dt + \frac{\sigma'(x_*)}{2\sqrt{\sigma(x_*)}} x_1 \ dW(t),\dots \label{eq:SDE3} \\
\end{align}

Equation~\eqref{eq:SDE1} is the deterministic equation and therefore
$x_*$ is an equilibrium of the deterministic equation.
Equations~\eqref{eq:SDE2} and~\eqref{eq:SDE3} are Ornstein--Uhlenbeck
processes. Using the well known solution to an Ornstein--Uhlenbeck
process and the It\^{o} isometry we can calculate the centered moments
of Equations~\eqref{eq:SDE2} and~\eqref{eq:SDE3} to obtain the first
and second order corrections to moments calculated in the linear noise
approximation:
\begin{align}
&\E[x] = x_* + \frac{f''(x_*)\sigma(x_*)}{4N(f'(x_*))^2}, \label{eq:SDEmoment1} \\
&\Var(x) = -\frac{\sigma(x_*)}{2Nf'(x_*)} + \frac{(f''(x_*))^2 \sigma^2(x_*)}{8N^2(f'(x_*))^4} + \frac{(\sigma'(x_*))^2}{16N^2(f'(x_*))^2}, \label{eq:SDEmoment2} \\
&\E[(x - \E[x])^3] = \frac{(f''(x_*))^3 \sigma(x_*)^3}{8N^3f'(x_*)^6} + \frac{3(\sigma'(x_*))^2 f''(x_*) \sigma(x_*)}{32 N^3 f'(x_*)^4}. \label{eq:SDEmoment3}
\end{align}

We can find any higher-order terms desired by continuing to solve the
system of SDEs iteratively.

Note that to leading order, these are the moments obtained from the van Kampen approximation.

\subsection{Laplace's method}\label{app:Laplace}

Equation~\eqref{eq:odetau} can be solved
using integrating factors to find (q.v.~\eqref{eq:stationary3}),
\begin{equation} \label{eq:tau}
\tau(x) = 2N \int_x^{x_+} e^{N\Phi(x')} \int_0^{x'} \frac{e^{-N\Phi(x'')}}{B(x'')} \ dx'' \ dx'.
\end{equation}
For large $N$, $e^{-N\Phi(x)}$ is sharply peaked around $x_-$ and
therefore $\int_0^{x'} e^{-N\Phi(x'')} \ dx''$ is nearly constant $x'
\in [x, x_+]$ for which $e^{N\Phi(x')}$ is well above zero. With this
assumption, equation~\eqref{eq:tau} can be written as
\begin{equation}
\tau(x) = 2N \int_x^{x_+} e^{N\Phi(x')} \ dx' \int_0^{x_+} \frac{e^{-N\Phi(x'')}}{B(x'')} \ dx''.
\end{equation}
The first integral is sharply peaked around $x_0$ and the second
integral is sharply peaked around $x_-$. Using a Gaussian
approximation (that is, replacing $\Phi(\cdot)$ with its second-order
Taylor expansion), and replacing the upper limits with any point
strictly above $x_0$, we obtain the MFPT from $x_-$ to $x_+$ to be

\subsection{WKB details}\label{app:WKB}

 the following $O(1)$ and $O(\frac{1}{N})$
equations
\begin{equation} \label{eq:order1}
0 = \Omega_+(x)\left(e^{\Psi'(x)}-1 \right) + \Omega_-(x)\left(e^{-\Psi'(x)} -1 \right),
\end{equation}
and
\begin{equation}\label{eq:order1N}
\begin{split}
  0 &= \Omega_+(x) e^{\Psi'(x)} \left( -\frac{k'(x)}{k(x)} - \frac{\Psi''(x)}{2} \right) + \Omega_-(x) e^{-\Psi'(x)} \left( \frac{k'(x)}{k(x)} - \frac{\Psi''(x)}{2} \right)\\&\quad - \Omega_+'(x) e^{\Psi'(x)} + \Omega_-'(x) e^{-\Psi'(x)}.
\end{split}
\end{equation}
The assumption that $\tau_-$ is large enough to not appear in these
equations is not yet justified, but we show below that this is a
self-consistent choice.  Solving~\eqref{eq:order1} for $\Psi(x)$ gives
two solutions:
\begin{equation}
\Psi_1(x) = const, \quad\Psi_2(x)=\int_0^x \ln\left( \frac{\Omega_-(x')}{\Omega_+(x')}\right) \ dx'.
\end{equation}
Substituting each into the next order equation~\eqref{eq:order1N} and
solving for $k(x)$ gives
\begin{equation}
k_1(x) = \frac{B}{\Omega_+(x) - \Omega_-(x)},\quad  k_2(x) =  \frac{A}{\sqrt{\Omega_+(x)\Omega_-(x)}}
\end{equation}



We now seek to approximate $\tau_-$ in the quasistationary
solution~\eqref{eq:defofPi}.  Let us start with the Ansatz that the
WKB approximation takes the form of the activation solution on $[0,
x_0]$ as the process tries to escape the basin of attraction around
$x_-$ and takes the form of the relaxation solution on $[x_0, 1]$ as
the process is quickly attracted to $x_+$.  The relaxation and
activation solutions need to satisfy the matching boundary conditions
at $x_0$, thus we match both to an appropriate inner solution at $x_0$
to find $A$ and $B$.  Begin by rewriting the diffusion approximation
as
\begin{equation} \label{eq:diffusionFlux}
\frac{\partial P(x,t)}{\partial t} = - \frac{\partial J(x,t)}{\partial x}
\end{equation}
where $$J(x, t) = (\Omega_+(x) - \Omega_-(x)) P(x, t) - \frac{1}{2N}
\frac{\partial}{\partial x} \left[ (\Omega_+(x) + \Omega_-(x) ) P(x,
  t) \right].$$

Substituting the quasistationary solution $\Pi(x, t) = \widehat\Pi(x)
e^{-t/\tau_-}$ into equation~\eqref{eq:diffusionFlux} and Taylor
expanding around $x_0$ to order $1/N$ (again assuming that $\tau_-\gg N$)
results in the constant flux through $x_0$,
\begin{equation} \label{eq:JconstExp}
J_0 = (x-x_0)(\Omega_+'(x_0) - \Omega_-'(x_0)) \widehat\Pi(x) - \frac{1}{N}(\Omega_+(x_0)) \widehat\Pi'(x)
\end{equation}
where we have used the fact that $\Omega_+(x_0) - \Omega_-(x_0) = 0$. It is straightforward to solve equation \eqref{eq:JconstExp} with an integrating factor to find the inner solution,
\begin{equation} \label{eq:inner}
\widehat\Pi_{inner}(x) = \frac{JN}{\Omega_+(x_0)} e^{\frac{(x-x_0)^2}{\sigma^2}} \int_x^\infty e^{-\frac{(y-x_0)^2}{\sigma^2}} \ dy
\end{equation}
where $$\sigma = \sqrt{\frac{2\Omega_+(x_0)}{N(\Omega_+'(x_0) -
    \Omega_-'(x_0))}}$$ determines the size of the boundary layer. In
order to match this inner solution with the activation and relaxation
solutions, we use the following asymptotics of the inner solution:
\begin{equation} \label{innerasym}
\widehat\Pi_{inner}(x) = \begin{cases} \dfrac{NJ_0\sigma^2}{2(x-x_0)\Omega_+(x_0)} & x-x_- \gg \sigma \\
\dfrac{J_0N\sigma \sqrt{\pi}}{\Omega_+(x_0)} e^{\frac{(x-x_0)^2}{\sigma^2}}  & x_0-x \gg \sigma \end{cases}
\end{equation}

Expanding the relaxation solution \eqref{eq:relaxation} around $x_0$
and matching to the inner solution for $x-x_0 \gg \sigma$ gives
$J_0=B$. Expanding the activation solution around $x_0$ and matching
to the inner solution for $x_0-x \gg \sigma$ gives

\[ J_0 = \frac{A\Omega_+(x_0)}{\sqrt{\Omega_+(x_0) \Omega_-(x_0)}} \sqrt{ \frac{|\Psi''(x_0)|}{2 \pi N}} e^{-N\Psi(x_0)} \]

Substituting the quasistationary solution $\Pi(x, t)$ into the diffusion approximation and integrating over $[0, x_0]$, we link $\tau_-$ and $J_0$
\begin{equation}
\tau_- = \frac{1}{J_0} \int_0^{x_0} \widehat\Pi(x) \ dx
\end{equation}

Using the activation solution on $[0, x_0]$ and a Gaussian approximation around $x_-$ results in
\begin{equation}
\tau_- = \frac{2\pi}{\Omega_+(x_-)} \frac{1}{\sqrt{|\Psi''(x_0)|\Psi''(x_-)}} e^{N(\Psi(x_0)-\Psi(x_-))}
\end{equation}

\end{document}